\documentclass[11pt]{amsart}
\usepackage{amssymb}
\usepackage{amsfonts}
\usepackage{amsthm}

\newenvironment{Proof}{\noindent{\begin{minipage}[t]{5.9 in}}{\it Proof.  }}
{\noindent{$\diamondsuit$}{\end{minipage}}}

\newtheorem{thm}{Theorem}[section]

\newtheorem{prop}{Proposition}[section]

\theoremstyle{remark}

\newtheorem*{ack}{Acknowledgments}

\newcommand{\spec}{\mathop{\rm Spec}\nolimits}
\newcommand{\dvs}{\mathop{\rm Div}\nolimits}

\newcommand{\hdeg}{\mathop{\widehat {\rm deg}}\nolimits}

\newcommand{\Q}{{\mathbb{Q}}}

\newcommand{\p}{{\mathbb{P}}}

\newcommand{\R}{{\mathbb{R}}}

\newcommand{\cA}{\mathcal A}
\newcommand{\cB}{\mathcal B}
\newcommand{\cX}{\mathcal X}
\newcommand{\cV}{\mathcal V}
\newcommand{\cU}{\mathcal U}

\newcommand{\cL}{\mathcal L}
\newcommand{\bL}{\overline {\cL}}
\newcommand{\gA}{{\mathcal A}_{\eta}}
\newcommand{\rA}{A_{\rho}}
\newcommand{\rL}{D_{\rho}}
\newcommand{\sx}{{\sigma}_x}

\newcommand{\bK}{{\overline K}}
\newcommand{\bH}{{\overline H}}
\newcommand{\bZ}{{\overline Z}}
\newcommand{\hg}{{h_{K(C)}^{\rm {geom}}}}
\newcommand{\chgl}{{{\hat h}_{({\rA}, {\rL})}^{\rm {geom}}}}
\newcommand{\ha}{{h_L^{\rm {arith}}}}
\newcommand{\hal}{{h_{(A, L)}^{\rm {arith}}}}
\newcommand{\hax}{{h_{(X, L)}^{\rm {arith}}}}
\newcommand{\chg}{{{\hat h}_{K(C)}^{\rm {geom}}}}
\newcommand{\cha}{{{\hat h}_L^{\rm {arith}}}}
\newcommand{\chat}{{{\hat h}_{(A_t, L_t)}^{\rm {arith}}}}
\newcommand{\hc}{\widehat c}

\newcommand{\maprlim}[1]{\smash{\mathop{\longrightarrow}\limits^{#1}}}

\begin{document}

\title{On the Specialization Theorem for Abelian Varieties} 
\author{Rania Wazir}
\email{wazir@dm.unito.it}
\address{Dipartimento di Matematica, Universit\`a degli Studi di Torino, Via Carlo Alberto 10, 10123 Torino, Italy.}
\subjclass{11G40, 14K15}
\keywords{Abelian fibrations, Mordell-Weil rank, Specialization Theorems, arithmetic height functions}

\maketitle

\bibliographystyle{amsalpha} 
\nocite{*}

\begin{abstract}
In this note, we apply Moriwaki's arithmetic height functions to obtain an analogue of Silverman's Specialization Theorem for families of Abelian varieties over $K$, where $K$ is any field finitely generated over ${\Q}$.
\end{abstract}

\par
Let $k$ be a number field and \mbox{${\pi}: {\cA} \rightarrow {\cB}$} be a proper flat morphism of smooth projective varieties over $k$ such that the generic fiber ${\gA}$ is an Abelian variety defined over $k({\cB})$.  For almost all absolutely irreducible divisors $D/k$ on ${\cB}$, \mbox{${\cA}_D := {\cA} {\times}_{\cB} D$} is also a flat family of Abelian varieties, and our goal is to compare the Mordell-Weil rank of ${\gA}$ with the Mordell-Weil rank of the generic fiber of ${\cA}_D$.  In fact, if we fix a projective embedding of ${\cA}$, ${\cB}$ into ${\p}^N$, and an integer $M > 0$, then we can show that for all but finitely-many divisors $D$ of degree less than $M$, the rank of the generic fiber of ${\cA}_D$ is at least the rank of ${\gA}$.  

\par
This result is a an amusing example of the height machine in action.  The main issue is to rephrase the problem in terms of Abelian varieties over function fields, and define the ``right'' height functions; the proofs then follow verbatim as in {\cite{js3}}.

\par
Assume now that we also have a proper, flat morphism \mbox{$ f: {\cB} \rightarrow {\cX}$} with generic fiber a smooth, irreducible curve $C$ defined over $K := k({\cX})$.  Composition gives a flat morphism \mbox{$ g := f \circ {\pi}: {\cA} \rightarrow {\cX} $} whose generic fiber is a smooth, irreducible variety $A$ defined over $K$, and by base extension we have a flat morphism \mbox{$ {\rho}: A \rightarrow C$}, whose generic fiber is ${\gA}$.  Observe that in this setting, divisors $D \in \dvs({\cB})$ such that $f(D) = {\cX}$ correspond to points on $C$.  The next proposition (based on a classical geometric argument, see {\cite[Proposition 5.1]{hp}}) shows that, up to replacing ${\cB}$ and ${\cA}$ by birationally equivalent varieties, we can always reduce to this situation.

\begin{prop}
Let ${\cV}$ be a smooth projective variety of dimension $n$ defined over $k$.  Then, there exist finitely many birationally equivalent varieties ${\nu}_i: {\tilde {\cV}}_i \rightarrow {\cV}$ and proper flat morphisms $f_i: {\tilde {\cV}}_i \rightarrow {\p}^{n-1}$ such that: 
\begin{enumerate}
\item
The generic fiber of $f_i$ is a smooth, irreducible curve $C_i$ defined over $k({\p}^{n-1})$.
\item
Let ${\cU}_i$ be an affine open dense subset of ${\cV}$ such that ${\nu}_i$ is an isomorphism on ${\cU}_i$, and let ${\cU} : = \cap \, {\cU}_i$.  Then for all divisors $D \in \dvs({\cV})$ such that $D \cap {\cU} \neq \emptyset$, there exists an $i$ such that $f_i(D') = {\p}^{n-1}$, where $D' := {\overline {{\nu}_i^{-1}(D \cap {\cU}_i)}}$.
\end{enumerate}
\label{geom}
\end{prop}

\begin{Proof}
Consider an embedding ${\cV} \hookrightarrow {\p}^N$, and take a linear projection inducing a finite morphism \mbox{$ {\phi}: {\cV} \rightarrow {\p}^n$}.  Let $p_0 \in {\p}^n$ be a point outside the ramification locus of ${\phi}$, let 
\mbox{${\mu}_0: {\tilde {\p}}_0 \rightarrow {\p}^n$} be the blow-up at $p_0$, and 
\mbox{${\mu}_0': {\tilde {\p}}_0 \rightarrow {\p}^{n-1}$} be the associated morphism. If ${\nu}_0:{\tilde {\cV}}_0 \rightarrow {\cV}$ is the blow-up of ${\cV}$ at the points ${\phi}^{-1}(p_0)$, then we have a morphism 
\mbox{$ {\tilde {\phi}}_0: {\tilde {\cV}}_0 \rightarrow {\tilde {\p}}_0$} which, when composed with ${\mu}_0'$, gives a fibration in curves 
\mbox{$f_0 := {\mu}_0' \circ {\tilde {\phi}}_0 : {\tilde {\cV}}_0 \rightarrow {\p}^{n-1}$}.  

\par
Next, let $D'$ be the proper transform of $D$ under the blow-up ${\mu}_0$.  Then \mbox{$f_0(D') \neq {\p}^{n-1}$} if and only if $D$ is a cone with vertex \mbox{$V_D \ni p_0$}.  Since the vertex of a cone is a linear subspace, the highest dimensional vertex is a hyperplane.  Therefore, if we choose $n+1$ points $p_0, \dots p_n$ not all lying in a hyperplane, we see that for all divisors $D$ in ${\p}^n$, we can find a morphism $f_i$ such that $f_i(D') = {\p}^{n-1}$.  Since ${\phi}$ is a finite morphism, the same holds for divisors in ${\cV}$.
\end{Proof}

\par
The geometric proposition above allows us to reduce our problem to the following situation:  Let $K$ be a field finitely generated over ${\Q}$, and let $C$ be a smooth projective curve defined over $K$.  Suppose $A$ is a smooth projective variety equipped with a proper flat morphism \mbox{${\rho}: A \rightarrow C$}, such that the generic fiber ${\rA}$ is an Abelian variety with Chow trace \mbox{$({\tau, B})$}.  Let \mbox{$ x \in C({\bK})$} be a point for which the fiber $A_x$ is nonsingular.  Then the specialization map 
$$
{\sigma}_x : A(C/K) \rightarrow A_x({\bK})
$$
is a homomorphism from the group of sections, to the group of points on the fiber. 

\begin{thm} Let ${\Gamma}$ be a finitely generated free subgroup of $A(C/K)$ which injects in $A(C/K)/{\tau}B(K)$.  Then the set 
$$
\left\{ x \in C({\bK}): {\sx} \; \textrm{is not injective on ${\Gamma}$} \right\}
$$
is a set of bounded height in $C({\bK})$.  In particular, if $d \geq 1$, then ${\sx}$ is injective for all but finitely many \mbox{$ x \in {\cup}_{[L:K] \leq d} C(L)$}.  Furthermore, if ${\rA}$ has trivial Chow trace, this shows that, excluding a finite number of points $x \in C$, the Mordell-Weil rank of the special fibers $A_x$ is at least that of the generic fiber ${\rA}$.
\label{thm1}
\end{thm}

\par
The proof of this theorem is based on being able to measure the variation of certain height functions in a family of Abelian varieties.  We briefly describe the necessary height functions below, before outlining the major steps in the proof.

\par
{\it Arithmetic Height.}  We consider the arithmetic height functions on $K$ introduced by Moriwaki {\cite{am1}}, and,  as much as possible, stick to Moriwaki's notation and terminology.  Let $Z$ be a normal projective arithmetic variety whose function field is $K$, and fix nef \mbox{$C^{\infty}$-hermitian} line bundles \mbox{$ {\bH}_1, {\bH}_2, ... , {\bH}_d$} on $Z$.  The collection \mbox{$(Z; {\bH}_1, {\bH}_2, ... , {\bH}_d)$} is called a polarization of $Z$, and will be denoted by ${\bZ}$.

\par
Now suppose that $X$ is a smooth projective variety over $K$, and $L$ a line bundle on $X$.  If ${\cX}$ is a projective arithmetic variety over $Z$ and ${\bL}$ a hermitian line bundle on ${\cX}$ with \mbox{${\cX}_K = X$} and \mbox{${\cL}_K = L$} then the pair $({\cX}, {\bL})$ is a model for $(X, L)$, and we can define a height function using Arakelov intersection theory, as follows:  \mbox{$h^{\bZ}_{({\cX}, {\bL})} : X({\bK}) \rightarrow {\R}$} via
$$
h^{\bZ}_{({\cX}, {\bL})} (P) := \frac{{\hdeg}\left({\hc}_1(\left.{\bL}\right|_{\Delta_P})\cdot
{\hc}_1(\left.f^*{\bH}_1\right|_{\Delta_P})\cdots
{\hc}_1(\left.f^*{\bH}_d\right|_{\Delta_P})\right)}{[K(P):K]},
$$
\noindent
where ${\Delta_P}$ is the Zariski closure of the point $P$ in 
\mbox{$ \spec ({\bK}) \maprlim{P} X \hookrightarrow {\cX} $}, and \mbox{$f: {\cX} \rightarrow Z$} is the canonical morphism.  Furthermore, if the polarization $Z$ is big (as defined in {\cite[Section 2]{am1}}), then $h^{\bZ}_{({\cX}, {\bL})}$ is an arithmetic height function, and denoted by $h^{\rm{arith}}_{({\cX}, {\bL})}$.  In fact, {\cite[Corollary 3.3.5]{am1}} shows that, up to a bounded function, the height $h^{\bZ}_{({\cX}, {\bL})}$ does not depend on the choice of model $({\cX}, {\bL})$ for $(X, L)$, and hence we will denote the arithmetic height function simply by ${\hax}$, or ${\ha}$ when $X$ is understood.

\begin{prop}[Height Machine]
Let $X$ be a smooth projective variety over $K$, and $L$ a line bundle on $X$.  Then the arithmetic height function ${\ha}$ satisfies the following properties:
\begin{enumerate}
\item
If $M$ is another line bundle on $X$, then \mbox{$h_{L{\otimes}M}^{\rm{arith}} = {\ha} + h_M^{\rm{arith}} + O(1)$}, while \mbox{$h_{L^{{\otimes}-1}}^{\rm{arith}} = -{\ha} + O(1)$}.
\item
Let $Bs(L)$ denote the base locus of $L$, and set \mbox{$SBs(L) := \cap_{n>0} Bs(L^{{\otimes}n})$}.  Then ${\ha}$ is bounded below on $(X{\backslash}SBs(L))({\bK})$, and in particular:
\begin{enumerate}
\item
If $L$ is ample, then ${\ha}$ is bounded below.
\item
If $L = {\mathcal O}_X$ then ${\ha} = O(1)$.
\end{enumerate}
\item
If $N$ is any number, and $R$ a positive integer, then the set
$$
\left\{ P \in X({\bK}) | {\ha}(P) \leq N \; \, \textrm{and} \; \: [K(P):K] \leq R \right\}
$$
\noindent
is finite.
\item
Let \mbox{$q: X \rightarrow Y$} be a morphism of smooth, projective varieties over $K$, and $M$ a line bundle on $Y$; then
$$
h_{(X, q^*(M))}^{\rm{arith}}(P) = h_{(Y, M)}^{\rm{arith}}(q(P)) \qquad \textrm{for all $P \in X({\bK})$.}
$$
\item
If $L$ is ample, and $M$ is any line bundle on $X$, then there is a constant $c$ such that
$$
{h_M^{\rm{arith}}} < c{\ha} + O(1).
$$
\item
Suppose $X$ is a smooth projective curve over $K$, and if $D$ is a divisor write $h_D^{\rm{arith}}$ for $h_{L(D)}^{\rm{arith}}$, where $L(D)$ is the line bundle associated to $D$. If $D$ is of degree $d>0$, and $E$ is a divisor of degree $e$, then
$$
\lim_{h_D^{\rm{arith}}(t) \rightarrow \infty} \;
\frac{h_E^{\rm{arith}}(t)}{h_D^{\rm{arith}}(t)} = \frac{e}{d} \, .
$$
\end{enumerate}
\end{prop}

\begin{Proof}
Properties (1)-(3) are {\cite[Proposition 3.3.7]{am1}}, while Property (4) is {\cite[Theorem 4.3]{am1}}.  To prove Property (5), we start with the definition of arithmetic height, and note that if $f: Y \rightarrow B$ is the canonical morphism on $Y$, then the canonical morphism on $X$ is given by $f' = f \circ q$: 
\begin{eqnarray*}
h_{(X, q^*(M))} (P)  & := & \frac{{\hdeg}\left( ({f'}^*) \left({\hc}_1({\bH}_1) \cdots {\hc}_1({\bH}_d) \right)\cdot
{\hc}_1(q^*({\bar M}))\cdot ({\Delta}_P, 0) \right)}{[K(P):K]} \\
& = & \frac{{\hdeg}\left( (q^* \circ f^*) \left({\hc}_1({\bH}_1) \cdots {\hc}_1({\bH}_d) \right)\cdot
q^*{\hc}_1({\bar M})\cdot ({\Delta}_P, 0) \right)}{[K(P):K]} \\
& = & \frac{ \deg({\Delta}_P \rightarrow {\Delta}_{q(P)} ){\hdeg}\left( (f^*) \left({\hc}_1({\bH}_1) \cdots {\hc}_1({\bH}_d) \right)\cdot
{\hc}_1({\bar M})\cdot ({\Delta}_{q(P)}, 0) \right)}{[K(P):K]} \\
& = & h_{(Y, M)} (q(P)),
\end{eqnarray*}
where the third equation follows from the Projection formula ({\cite[Proposition 1.3]{am1}}), and the final step is again by definition.

\par
To prove Property (6), we note that there is some constant $a$ such that \mbox{$L^{{\otimes}a}{\otimes}M^{{\otimes}{-1}}$} is ample.  Hence, by Property (3.a) we have 
\mbox{$h_{L^{{\otimes}a}{\otimes}M^{{\otimes}-1}} > O(1)$}.  The result now follows by Properties (1) and (2).

\par
Finally, the proof of Property (7) is very similar to the proof for heights over number fields (see for example {\cite[Chapter 4, Corollary 3.5]{sl1}}).  
\end{Proof}

\par
If $A/K$ is an Abelian variety, then we can also define a canonical height ${\cha}$, which has the following properties:

\begin{prop}[Canonical Height Machine]
\begin{enumerate}
\item
If $M$ is another line bundle on $A$, then 
\mbox{${\hat h}_{L{\otimes}M}^{\rm{arith}} = {\cha} + {\hat h}_M^{\rm{arith}}$}, while 
\mbox{${\hat h}_{L^{{\otimes}-1}}^{\rm{arith}} = -{\cha} $}.
\end{enumerate}
\par
If we assume furthermore that $L$ is ample, then we have:
\begin{enumerate}
\item[(3)]
There is a constant $c$, such that if $M$ is another line bundle on $A$, then 
\mbox{$ {\hat h}_M^{\rm{arith}} \leq c{\cha}$}.
\item[(4)]
${\cha}(x) \geq 0$ for all $x \in A({\bK})$, and equality holds if and only if $x$ is a torsion point.
\end{enumerate}
\end{prop}

\begin{Proof}
This is {\cite[Propositions 3.4.1 and 3.4.2]{am1}}.  He proves these under the condition that $L$ is a symmetric line bundle, but the proofs hold verbatim for $L$ anti-symmetric, hence for all line bundles $L$. 
\end{Proof}

\par
{\it Geometric Height.}  Since $K$ is a global field, and $C/K$ is a smooth projective curve, $K(C)$ is also a global field, and we can define the geometric height ${\hg}$ so that, for any \mbox{$ x \in K(C) $}, \mbox{${\hg}([x,1])$} is the degree of the morphism \mbox{$ [x,1]:C \rightarrow {\p}^1 $}.  If $A$ is an Abelian variety defined over $K$, then we also have a canonical height ${\chg}$, which is non-degenerate on $A({\bK})/\left({\tau}B(K) + A_{{\rm {tors}}} \right)$ {\cite[Chapter 6, Theorem 5.4]{sl1}}.  Both ${\hg}$ and ${\chg}$ satisfy the usual properties of heights/canonical heights, see {\cite[Chapter 3]{sl1}}.

\par
Equipped with these height functions, the proof of the following theorems then proceeds exactly as in {\cite{js3}}, and leads to a proof of Theorem {\ref{thm1}}:

\par
Fix a line bundle $L$ on $A$, and let ${\hal}$ be the associated arithmetic height function.  For each \mbox{$t \in C({\bK})$} let $L_t$ be its restriction to $A_t$, and ${\rL}$ be its restriction to the generic fiber ${\rA}$.  Let $C^0 \subset C$ be an affine open subset such that for all \mbox{$t \in C^0({\bK})$}, the fiber $A_t$ is an Abelian variety.  Then there is a geometric canonical height 
$$
{\chgl}: {\rA}(K(C)) \rightarrow {\R},
$$
and for each \mbox{$t \in C^0({\bK})$}, an arithmetic canonical height
$$
{\chat}: A_t({\bK}) \rightarrow {\R}.
$$

\par
Set \mbox{$U := {\pi}^{-1}(C^0)$}.  The canonical heights ${\hat h}_{(A_t, L_t)}^{\rm{arith}}$ on the good fibers $A_t$ can be fitted together to give a "canonical height" on $U({\bK})$, \mbox{${\hat h}_L: U({\bK}) \rightarrow {\R}$}.  If we fix also a line bundle $M$ on $C$, then 

\begin{thm}
There is a constant $c$, depending on $M$, $L$, and the family $A \rightarrow C$, such that
$$
\left|{\hat h}_L(P) - {\ha}(P) \right| < ch_M^{\rm{arith}}(P) + O(1) \qquad 
\textrm{for all $P \in U({\bK})$.}
$$
(The $O(1)$ depends on the choice of heights ${\ha}$ and $h_M^{\rm{arith}}$ but not on $P$.)
\label{thm2}
\end{thm}

\par
Consider now the case where $C$ is a smooth projective curve over $K$, and fix an arithmetic height on $C$ as follows: Let $D$ be a divisor on $C$ with $\deg D > 0$, and define
$$
h_C^{\rm{arith}} := \frac{1}{\deg D} h_{(C, D)}^{\rm{arith}}.
$$
\begin{thm}
With notation as above, fix a section $P \in A(C)$.  Then
$$
\lim_{{t \in C^0({\bK})} \atop {h_C^{\rm{arith}}(t) \rightarrow \infty}} \frac{{\chat}(P_t)}{h_C^{\rm{arith}}(t)} = {\chgl}(P_{\rho}).
$$
\label{thm3}
\end{thm}
Note that by Property (7) of arithmetic heights, this result does not depend on the choice of arithmetic height on $C$.

\par
Finally, we observe that in the case that $A$ is an elliptic surface, another crank of the height machine yields the following sharpened estimate \'a la Tate {\cite{jt3}} for Theorem {\ref{thm3}}:

\begin{thm}
Assume that \mbox{${\pi}: A \rightarrow C$} is an elliptic surface, and assume that 
\mbox{$P \in A(C/K)$} is a section.  Then
$$
{\chat}(P_t) = {\chgl}(P_{\rho})h_C^{\rm{arith}}(t) + O \left( \sqrt{h_C^{\rm{arith}}(t)} + 1 \right)
$$
\label{thm4}
\end{thm}

\par
Tate's proof requires the existence of a good compactification of the Neron model, hence is not known to apply to families of higher dimensional Abelian varieties.  However, in the case of number fields, Call{\cite[Theorem I]{gc}}, using a different approach, generalized the result to the case of a one-parameter family of Abelian varieties.  It remains an interesting question whether an analogous estimate exists for a family of Abelian varieties over $K$.

\begin{ack}
The author would like to thank Alberto Albano, Joe Silverman, and Siman Wong for their useful advice and suggestions.
\end{ack}

\bibliography{spec} 
 
\end{document}